\documentclass [10pt,a4paper] {article}

\usepackage[french]{babel}

\usepackage{amssymb,amsmath}

\usepackage [all] {xy}

\usepackage [dvips] {graphicx}

\author{St\'ephane \textsc{Dugowson}\,\footnote{St\'ephane \textsc{Dugowson}, dugowson@ext.jussieu.fr.}}
\title {Representation of {f}inite connectivity spaces}
\date{\today} 

\newtheorem{df}{D\'efinition} 
 
\newtheorem{prop} {Proposition} 
 
\newtheorem{thm} [prop] {Th\'eor\`eme} 

\begin{document}

\maketitle

\mbox{}

 \textsc{Abstract} --- After recalling the definition of connectivity spaces and some of their main properties, a way is proposed to represent finite connectivity spaces by directed simple graphs. Then a connectivity structure is associated to each tame link. It is showed that all spaces of a certain class (the iterated Brunnian ones) admit representations by links. Finally, I conjecture that every finite connectivity space is representable by a link.

\mbox{}

 \textsc{R\'esum\'e} ---  Repr\'esentation des espaces connectifs finis --- Apr\`es avoir rappel\'e la d\'efinition des espaces connectifs et certaines de leurs principales propri\'et\'es, nous proposons une fa\c con de repr\'esenter les espaces connectifs finis par des graphes simples orient\'es, puis nous associons \`a tout entrelacs une structure connective. Nous montrons que tout espace d'une certaine classe (les espaces brunniens it\'er\'es) admet une repr\'esentation par entrelacs, et nous conjecturons finalement que tout espace connectif fini est repr\'esentable par entrelacs.

\mbox{}


\mbox{}

 \textsc{Mots cl\'es} --- Espace connectif. Graphe. Entrelacs. Borrom\'een. Brunnien. Conjecture.
 
\mbox{}


\mbox{}

\begin{tabular}{p{3.7cm}p{7.2cm}}
{} & \quad\emph{Il est remarquable qu'une figure aussi simple que celle du n\oe ud borrom\'een n'ait pas servi de d\'epart \`a -- une topologie}
\end{tabular}

\begin{flushright}
Jacques Lacan~\cite{Lacan:1973}
\end {flushright}

\section *{Introduction}

\`A notre connaissance, la notion d'espace connectif a d'abord \'et\'e introduite en 1981 par B\"orger~\cite{Borger:1981, Borger:1983}, dans le cadre de la th\'eorie des cat\'egories. Elle a \'et\'e red\'ecouverte, de fa\c con ind\'ependante, en 1988 par Matheron et Serra~\cite{MatheronSerra:1988k} dans le cadre de leurs recherches en morphologie math\'ematique, puis par nous-m\^eme en 2003~\cite{Dugowson:2003, Dugowson:2007c}, \`a partir d'une r\'eflexion sur la nature topologique du jeu de go.

\begin{figure} 
\begin{center}
\includegraphics [scale=0.2]{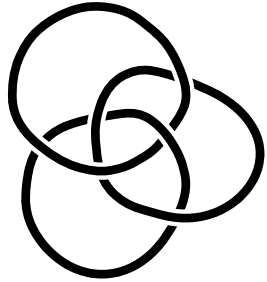}
\caption{Entrelacs borrom\'een}
\label{borro}
\end{center}
\end{figure}

Or, l'une des structures connectives parmi les plus simples, si elle n'est pas celle d'un espace topologique ni celle d'un graphe, pr\'esente une analogie \'evidente avec l'entrelacs borrom\'een (figure~\ref{borro}). Dans cet article, nous nous proposons de donner une forme math\'ematique pr\'ecise \`a cette analogie en associant \`a tout entrelacs un espace connectif fini (section~\ref{Representation par entrelacs}).  Nous montrons alors que toutes les structures connectives que nous appelons «~brunniennes it\'er\'ees~» sont associ\'ees \`a des entrelacs, et nous conjecturons finalement (section~\ref{Une conjecture}) que toute structure connective finie est celle d'un entrelacs.

Auparavant, afin de disposer d'outils permettant une description aussi simple que possible des espaces connectifs, nous introduisons dans la section~\ref{Representation par graphes simples orientes} les notions de \emph{partie connexe irr\'eductible} (section~\ref{Espaces connectifs irreductibles}), de \emph{graphe g\'en\'erique} (section~\ref{Graphe generique d un espace connectif}), et de \emph{familles connectivement libres} de parties d'un ensemble (section~\ref{Familles connectivement libres}). Nous signalons aussi certains r\'esultats quant au d\'enombrement des espaces connectifs finis ayant au plus six points, notamment ceux obtenus en 2002 par Wim van Dam~\cite{vanDam:2002a,vanDam:2002b} qui donne en particulier lieu \`a une suite de nombres premiers dont les cinq premiers termes --- les seuls connus \`a ce jour --- pr\'esentent une progression tr\`es rapide~\cite{Dugowson:2007suite}.

La section~\ref{Rappels sur les espaces connectifs} rappelle la d\'efinition de la cat\'egorie des espaces connectifs et ses principales propri\'et\'es.

L'article est illustr\'e par plusieurs figures. Toutes les repr\'esentations d'entrelacs ont \'et\'e produites \`a l'aide du logiciel \emph{KnotPlot} de Robert G. Scharein~\cite{Scharein:knotplot}.

\mbox{}

\noindent \emph{Notations utilis\'ees}.  Pour tout ensemble $X$, on note $\mathcal{P}(X)$ l'ensemble des parties de $X$, et $card(X)$ le cardinal de $X$. On note $X_n$ l'ensemble $\{1,\dots, n\}$. Etant donn\'e un ensemble $\mathcal{K}$ de parties d'un ensemble $X$, on notera $\mathcal{K}^*$ l'ensemble de ces parties ayant au moins deux \'el\'ements : $\mathcal{K}^*=\{K\in\mathcal{K}, card(K)>1\}$. Pour tout ensemble $X$, on pose $\mathcal{T}(X) = {\mathcal{P}(X)}^*$. Les symboles d'inclusion $\subset$ et $\supset$ sont pris au sens large.

\section {Rappels sur les espaces connectifs}\label{Rappels sur les espaces connectifs}

\begin {df} [Espaces connectifs]  Un \emph{espace connectif} est un  couple $\mathbf{X}=(X,\mathcal{K})$ form\'e d'un ensemble $X$ et d'un ensemble $\mathcal{K}$ de parties non vides de $X$ tel que 
\begin{displaymath}
\forall \mathcal{I}\in \mathcal{P}(\mathcal{K}), \bigcap_{K\in\mathcal{I}}K\ne\emptyset\Longrightarrow \bigcup_{K\in\mathcal{I}}K\in\mathcal{K}.
\end{displaymath}
\noindent L'ensemble $X$ est le \emph{support} de $\mathbf{X}$, l'ensemble $\mathcal{K}$ est une \emph{structure connective} sur $X$ et ses \'el\'ements sont les \emph{parties connexes} ou \emph{parties connect\'ees} de l'espace connectif $\mathbf{X}$.
Nous dirons d'un espace connectif qu'il est \emph{int\`egre} si tout  singleton est connect\'e.
Un \emph{morphisme connectif}, ou \emph{application connective}, d'un espace connectif $(X,\mathcal{K})$ vers un autre $(Y,\mathcal{L})$ est une application $f:X\to Y$ telle que :
\begin{displaymath}
\forall K\in \mathcal{K}, f(K)\in \mathcal{L}.
\end{displaymath}
\end{df}

Dans la suite de cet article, tous les espaces connectifs consid\'er\'es seront suppos\'es int\`egres. Pour un tel espace $(X,\mathcal{K})$, $X$ est enti\`erement d\'etermin\'e par la donn\'ee de $\mathcal{K}$, puisque $X=\bigcup_{K\in\mathcal{K}}K$, tandis que $\mathcal{K}$ est d\'etermin\'e par la donn\'ee de $X$ et de $\mathcal{K}^*$, puisque $\mathcal{K}=\bigcup_{x\in X}\{\{x\}\}\cup \mathcal{K}^*$.

\begin{thm} [Composantes connexes] Pour tout espace connectif (int\`egre) non vide, les composantes connexes, c'est-\`a-dire les parties connexes maximales, forment une partition.
\end{thm} 

 \noindent\emph{D\'emonstration}. L'espace \'etant int\`egre, tout point appartient \`a une partie connexe maximale, \`a savoir l'union des connexes qui contiennent ce point. De plus, les parties connexes maximales sont trivialement deux \`a deux disjointes.
\begin {flushright}
$\square$
\end{flushright}
\pagebreak[3]

\begin{prop} Ordonn\'e par l'inclusion, l'ensemble $\mathbf{K}(X)$ des structures connectives (int\`egres) sur un ensemble donn\'e de points $X$ constitue un treillis complet.
\end{prop} 

La cat\'egorie des espaces connectifs int\`egres constitue en fait ce que nous avons appel\'e une «~cat\'egorie \`a treillis de structures~”~\cite{Dugowson:2007c}, de sortes qu'elle admet toutes limites et colimites. C'est une cat\'egorie topologique, non cart\'esienne ferm\'ee~\cite{Borger:1983} mais mono\"\i dale ferm\'ee~\cite{Dugowson:2003, Dugowson:2007b, Dugowson:2007c}.

De la propositon pr\'ec\'edente, il d\'ecoule que tout ensemble $\mathcal{A}$ de parties de $X$ est contenue dans une plus petite (plus \emph{fine}) structure connective (int\`egre) $[\mathcal{A}]$ sur $X$. Le th\'eor\`eme suivant  (dans lequel $\omega_0$ d\'esigne le plus petit ordinal infini) conduit \`a appeler $[\mathcal{A}]$ la structure connective (int\`egre) \emph{engendr\'ee} par $A$.

\begin{thm} [Engendrement des structures connectives] Soit $X$ un ensemble, et $\mathcal{A}_0$ un ensemble de parties de $X$. Il existe un ordinal $\alpha_0$ inf\'erieur ou \'egal \`a $\omega_0+1$ tel que 
\begin{displaymath}
[\mathcal{A}_0]=\Psi^{\alpha_0}(\mathcal{A}_0),
\end{displaymath}
o\`u, pour tout ordinal $\alpha\ge1$, $\Psi^\alpha$ est l'application de $\mathcal{P}(\mathcal{P}(X))$ dans lui-m\^eme d\'efinie par induction pour tout $\mathcal{A}\in\mathcal{P}(\mathcal{P}(X))$ par
\begin{displaymath}
\begin{array}{l}
\bullet  \textrm{ } \Psi^1(\mathcal{A})=\Psi(\mathcal{A})=\{A\in\mathcal{P}(X), \exists\mathcal{L}\subseteq\mathcal{A}, \bigcap_{L\in\mathcal{L}}L\ne\emptyset\textrm{ et } A=\bigcup_{L\in\mathcal{L}}L \}.\\
\bullet \textrm{ Si }\alpha=\beta+1, \Psi^\alpha=\Psi\circ\Psi^\beta,\\
\bullet \textrm{ Si }\alpha\textrm{ est un ordinal limite,} \Psi^\alpha(\mathcal{A})=\bigcup_{\beta<\alpha}\Psi^\beta(\mathcal{A}).\\
\end{array}
\end{displaymath}
\end{thm} 

 \noindent\emph{D\'emonstration}. Voir ~\cite{Dugowson:2003}, p.\,11 - 13.
\begin {flushright}
$\square$
\end{flushright}
\pagebreak[3]

\begin{df} L'\emph{union brunnienne} $\biguplus_i\mathbf{Y}_i$ d'une famille non vide $\mathbf{Y}_i=(Y_i,\mathcal{K}_i)$ d'espaces connectifs non vides est l'espace connectif $\mathbf{X}=(X,\mathcal{K})$ de support $X=\bigsqcup_i Y_i$ (union disjointe des supports $Y_i$), et dont la structure connective est obtenue en ajoutant aux connexes des $\mathbf{Y}_i$ la partie pleine du nouvel espace : $\mathcal{K}= \bigsqcup_i \mathcal{K}_i \cup \{X\}$.
\end{df}

\begin{df} Pour tout entier $n\ge 1$, on appelle \emph{espace brunnien \`a $n$ points} l'union brunnienne de $n$ espaces r\'eduits \`a un point. On appelle en particulier \emph{espace borrom\'een} l'espace brunnien \`a trois points. On appelle \emph{structure brunnienne} (resp. \emph{borrom\'eenne}) la structure connective de ces espaces.
\end{df}

Ainsi, l'espace brunnien \`a $n$ points est l'espace dont le support $X$ comporte $n$ points et dont la structure connective est la structure brunnienne $\mathcal{K}$ caract\'eris\'ee par $\mathcal{K}^*=\{X\}$.

 \emph{Espaces brunniens it\'er\'es}. On d\'efinit par r\'ecurrence les espaces brunniens d'ordre $r\in\mathbf{N}$ : un espace brunnien d'ordre $0$ est un espace r\'eduit \`a un point; pour tout $r$, un espace brunnien d'ordre $r+1$ est l'union brunnienne d'un espace brunnien d'ordre $r$ et d'un ou plusieurs autres espaces brunniens d'ordres inf\'erieurs ou \'egaux \`a $r$.

\section {Repr\'esentation par graphes simples orient\'es}\label{Representation par graphes simples orientes}

\subsection {Espaces connectifs irr\'eductibles}\label {Espaces connectifs irreductibles}

\begin {df} Une partie connexe $R$ d'un espace connectif est dite \emph{r\'eductible} si elle est l'union de deux parties propres connexes d'intersection non vide :
\begin{displaymath}
R\textrm{ r\'eductible }\Leftrightarrow\exists(K_1,K_2)\in\mathcal{K}^2,(\forall i\in\{1,2\},K_i\subsetneqq R) \textrm{ et }  (K_1\cap K_2 \ne \emptyset).
\end{displaymath}
Une partie connexe non r\'eductible est dite \emph{irr\'eductible}. Un espace connectif est dit irr\'eductible si son support est un connexe irr\'eductible.
\end {df}

\begin {prop} Un espace connectif avec au moins deux points est irr\'eductible si et seulement s'il est l'union brunnienne de certains de ses sous-espaces propres.
\end {prop}

 \noindent\emph{D\'emonstration}. Toute union brunnienne est clairement irr\'eductible. R\'eciproquement, soit $(X,\mathcal{K})$ un espace connectif irr\'eductible tel que $card(X)\ge 2$. Alors $(X,\mathcal{K}\smallsetminus\{X\})$ est encore un espace connectif, dont les composantes connexes ont pour union brunnienne $(X,\mathcal{K})$.
\begin {flushright}
$\square$
\end{flushright}
\pagebreak[3]

\subsection {Graphe g\'en\'erique d'un espace connectif}\label {Graphe generique d un espace connectif}

A tout espace connectif int\`egre $\mathbf{X}=(X,\mathcal{K})$, on associe un graphe simple orient\'e $\mathbf{G}_\mathcal{K}=(G_\mathcal{K},\mathcal{A}_\mathcal{K})$ de la fa\c con suivante. L'ensemble $G_\mathcal{K}$ des sommets de $\mathbf{G}_\mathcal{K}$ est constitu\'e des parties connexes irr\'eductibles de $\mathbf{X}$, y compris les singletons. Les \'el\'ements de $G_\mathcal{K}$ seront \'egalement appel\'es \emph{points g\'en\'eriques} de $\mathbf{X}$. L'ensemble $\mathcal{A}_\mathcal{K}$ des ar\^etes orient\'ees de $\mathbf{G}_\mathcal{K}$ est form\'e des couples de points g\'en\'eriques distincts $(a,b)$ tels que $a\supset b$ et qu'il n'existe pas de point g\'en\'erique $c$ distinct de $a$ et de $b$ tel que $a\supset c\supset b$.

\begin{figure} 
\begin{center}
$$
\xymatrix{
  & x \ar[ld]  \ar[d] \ar[rd]&   \\
a  & b  & c \\
}
$$
\caption{Graphe g\'en\'erique de l'espace connectif borrom\'een}
\label{gg}
\end{center}
\end{figure}
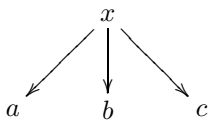

\begin{df} Le graphe simple orient\'e $\mathbf{G}_\mathcal{K}$ est le \emph{graphe g\'en\'erique} de l'espace connectif $\mathbf{X}$.
\end{df}

\begin{prop} \label{carKG} Tout espace connectif (int\`egre) \emph{fini} $(X,\mathcal{K})$ est caract\'eris\'e par son graphe g\'en\'erique $(G_\mathcal{K},A_\mathcal{K})$.
\end{prop}

\noindent \emph{Remarque}. L'hypoth\`ese de finitude est essentielle. Par exemple, la droite num\'erique usuelle n'admet que les singletons pour parties connexes irr\'eductibles, et sa structure connective n'en est pas moins distincte de la structure connective discr\`ete (\emph{i.e.} la plus fine, pour laquelle les seuls connexes sont les singletons).

\noindent\emph{D\'emonstration}. Soit $(X,\mathcal{K})$ un espace connectif int\`egre fini. On a $X=\bigcup_{u\in G_\mathcal{K}}u$. Par ailleurs, on v\'erifie par r\'ecurrence sur le cardinal (fini) des parties consid\'er\'ees que toute partie connexe de $X$ est l'union des parties connexes irr\'eductibles qu'elle contient, de sorte que $\mathcal{K}=[G_\mathcal{K}]$.
\begin {flushright}
$\square$
\end{flushright}
\pagebreak[3]

\noindent\emph{Exemples}. 1. L'espace borrom\'een \`a trois points $a, b, c$ admet quatre points g\'en\'eriques, trois qui s'identifient aux points de l'espace, et un quatri\`eme, disons $x$, qui s'identifie \`a l'espace entier. Son graphe g\'en\'erique est donc repr\'esent\'e par la figure\,\ref{gg}.

2. L'espace connectif d\'efini par
\begin{displaymath}
X=\{a,b,c,d\}\textrm{ et }\mathcal{K}^*=\{\{a,b\},\{b,c\},\{a,b,c\}, X\}, 
\end{displaymath}
a un graphe g\'en\'erique $\mathbf{G}_\mathcal{K}$ qui comporte, outre les quatre singletons et le point g\'en\'erique $x$ de l'espace, deux points g\'en\'eriques, disons $u$ et $v$, s'identifiant respectivement aux connexes irr\'eductibles $\{a,b\}$ et $\{b,c\}$. Notons que $b$, $u$, $v$ et $x$ forment un cycle du graphe consid\'er\'e, qui n'est donc pas un arbre.

\begin{prop} Un espace connectif fini est un espace brunnien it\'er\'e si et seulement si son graphe g\'en\'erique est un arbre; de plus l'ordre du premier co\"\i ncide alors avec la hauteur du second.
\end{prop}

\noindent \emph{D\'emonstration}. Le sens direct se v\'erifie par r\'ecurrence sur l'ordre des espaces brunniens it\'er\'es, la r\'eciproque par r\'ecurrence sur la hauteur des arbres.
\begin {flushright}
$\square$
\end{flushright}
\pagebreak[3]

\noindent \emph{Remarque}. Si le graphe g\'en\'erique d'un espace connectif est un arbre, tout n\oe ud de celui-ci a au moins deux fils.

La notion d'ordre d'un espace brunnien it\'er\'e s'\'etend \`a tout espace connectif fini $\mathbf{X}$. Pour cela, on commence par  d\'efinir par r\'ecurrence l'\emph{ordre} des points g\'en\'eriques de $\mathbf{X}$ : les points d'ordre $0$ sont ceux qui ne constituent l'origine d'aucune ar\^ete du graphe g\'en\'erique, autrement dit ce sont les singletons de $X$; pour tout $n$, un point g\'en\'erique $u$ est d'ordre $n+1$ si et seulement si $n$ est l'ordre maximal des extr\'emit\'es $v$ des arcs ayant ce point pour origine. 

\begin{df} L'\emph{ordre} d'un espace connectif fini est l'ordre maximal de ses points g\'en\'eriques.
\end{df}

\subsection {Familles connectivement libres}\label {Familles connectivement libres}

\begin{df} On dit qu'un ensemble de parties $\mathcal{L}\in\mathcal{P}(\mathcal{T}(X))$ est \emph{(connectivement) libre} si pour tout $K\in\mathcal{L}$, $K$ est un connexe irr\'eductible de l'espace connectif $(X,[\mathcal{L}])$. 
\end{df}

Cette d\'efinition, qui s'\'etend imm\'ediatement aux \emph{familles} de parties, signifie intuitivement que dans une famille libre la connexit\'e de certaines parties n'entra\^\i ne pas celle des autres. Dans la suite, on note $\mathbf{F}(X) = \{\mathcal{L}\in\mathcal{P}(\mathcal{T}(X)), \mathcal{L}\textrm{ est libre}\}$. Les propositions suivantes d\'ecoulent imm\'ediatement des consid\'erations pr\'ec\'edentes.

\begin{prop} Si $\mathcal{L}\in\mathbf{F}(X)$ et si $A\in\mathcal{T}(X)\smallsetminus[\mathcal{L}]$, alors $\mathcal{L}\cup\{A\}\in\mathbf{F}(X)$.
\end{prop}

\begin{prop} Si $X$ est un ensemble fini, l'application $\mathbf{K}(X)\to\mathbf{F}(X)$ d\'efinie par   $\mathcal{K}\mapsto {G_\mathcal{K}}^*$ est une bijection, de r\'eciproque $\mathcal{L}\mapsto[\mathcal{L}]$.
\end{prop}

\subsection {Enum\'eration des espaces connectifs finis}\label {Enumeration des espaces connectifs finis}

On se propose maintenant de d\'ecrire rapidement un proc\'ed\'e d'\'enum\'eration  de l'ensemble des structures connectives (int\`egres) de support l'ensemble fini  $X_n=\{1,\dots,n\}$. D'apr\`es la proposition pr\'ec\'edente, cette \'enum\'eration est \'equivalente \`a celle de l'ensemble $\mathbf{F}(X_n)$. Commen\c cons par munir l'ensemble fini $\mathcal{T}(X_n)$ d'une relation d'ordre totale, par exemple la relation $\preceq$ d\'efinie en posant, pour tout $L\in \mathcal{T}(X_n)$, $\delta(L)=\sum_{k\in L} {2^k}$ et, pour tout couple $(L_1,L_2)$ de parties de $X_n$ \`a plus de deux \'el\'ements, $L_1\preceq L_2 \Leftrightarrow \delta(L_1)\leq \delta(L_2)$.

On munit alors l'ensemble $\mathbf{F}(X_n)$ d'une relation d'ordre partielle $\preccurlyeq$ en posant, pour tout couple $(\mathcal{L},\mathcal{K})$ d'ensembles libres de parties \`a au moins deux \'el\'ements de $X_n$,
\begin{displaymath}
\mathcal{L}\preccurlyeq\mathcal{K}\Longleftrightarrow
\left\{
\begin{array}{l} 
\mathcal{L}\subset \mathcal{K}\\
\forall K\in \mathcal{K}\smallsetminus\mathcal{L}, \forall L\in \mathcal{L}, L \preceq K
\end{array}
\right.
\end{displaymath}

Lorsque la relation $\mathcal{L}\preccurlyeq\mathcal{K}$ est satisfaite, nous disons que $\mathcal{K}$ \emph{compl\`ete} $\mathcal{L}$ \emph{\`a droite}. On consid\`ere alors l'application $\Phi : \mathbf{F}(X_n) \to \mathcal{P}(\mathbf{F}(X_n))$ qui \`a tout $\mathcal{L}\in\mathbf{F}(X_n)$ associe l'ensemble des familles libres qui compl\`etent $\mathcal{L}$ \`a droite : $\Phi(\mathcal{L})=\{\mathcal{K}\in\mathbf{F}(X_n), \mathcal{L}\preccurlyeq \mathcal{K}\}$. Posant, pour tout $\mathcal{L}\in\mathbf{F}(X_n)$, 
\begin{displaymath}
\sigma(\mathcal{L})=\{A\in\mathcal{T}(X_n), \forall L\in\mathcal{L}, L\preceq A\}\smallsetminus[\mathcal{L}],
\end{displaymath}
on a
\begin{displaymath}
\Phi(\mathcal{L})=\{\mathcal{L}\}\cup\bigcup_{A\in\sigma(\mathcal{L})}\Phi(\mathcal{L}\cup \{A\}).
\end{displaymath}

En particulier, lorsque $\sigma(\mathcal{L})=\emptyset$, on a $\Phi(\mathcal{L})=\{\mathcal{L}\}$. Ceci fournit alors un proc\'ed\'e r\'ecursif d'\'enum\'eration des ensembles $\Phi(\mathcal{L})$, et en particulier de $\Phi(\emptyset)=\mathbf{F}(X_n)$, donc de $\mathbf{K}(X_n)$.

En utilisant ces id\'ees, nous avons impl\'ement\'e un programme fond\'e sur une proc\'edure r\'ecursive, et nous avons notamment obtenu les r\'esultats suivants. Notant $s_n = card (\mathbf{K}(X_n)) = card (\mathbf{F}(X_n))$ le nombre de structures connectives distinctes sur un ensemble \`a $n$ \'el\'ements \'etiquet\'es, $c_n= card(\{\mathcal{K}\in\mathbf{F}(X_n), X_n\in[\mathcal{K}]\})$ le nombre de ces structures pour lesquelles l'espace entier est lui-m\^eme connexe, $f_n=\max_{\mathcal{L}\in \mathbf{F}(X_n)} card(\mathcal{L})$ le nombre maximal de parties connexes irr\'eductibles non r\'eduites \`a un point pour un support \`a $n$ points, on a

$s_1=1, s_2=2, s_3=12, s_4=420, s_5=254076$,

$c_1=1, c_2=1, c_3=8, c_4=378, c_5=252 000$,

$f_1=0, f_2=1, f_3=3, f_4=6, f_5=13$.

\mbox{}

Remarquons que, pour $n\ge 2$, le nombre $k_n= card(\{\mathcal{K}\in\mathbf{F}(X_n), X_n\in\mathcal{K}\})$ des structures connectives sur $X_n$ pour lesquelles l'espace entier est un connexe irr\'eductible v\'erifie $k_n=s_n - c_n$, puisque l'application $\mathcal{K}\mapsto\mathcal{K}\cup\{X_n\}$  r\'ealise une bijection entre les structures non connexes et les structures irr\'eductiblement connexes.

Le nombre $f_5=13$ est par exemple illustr\'e par la structure connective sur $X_5$ d\'efinie par $\mathcal{K}=$ $[\{1,2\},$ $\{1,3\},$ $\{2,3\},$ $\{1,2,4\},$ $\{1,3,4\},$ $\{2,3,4\},$ $\{1,2,5\},$ $\{1,3,5\},$ $\{2,3,5\},$ $\{4,5\},$ $\{1,4,5\},$ $\{2,4,5\},$ $\{3,4,5\}]$.

Le site \emph{On-Line Encyclopedia of Integer Sequences} de Neil J. A. Sloane contient en outre la valeur de $s_6$, calcul\'ee en 2002 par Wim van Dam~\cite{vanDam:2002a,vanDam:2002b} : $s_6 = 17199454920$. En notant $p_n$ le plus grand facteur premier de $s_n$, on remarque que l'on obtient une suite de nombres premiers dont les premiers termes croissent rapidement~\cite{Dugowson:2007suite}~: $p_2=2$, $p_3=3$, $p_4= 7$, $p_5 = 683$, $p_6=143328791$.

On doit \'egalement \`a Wim van Dam~\cite{vanDam:2002b} le calcul des premiers termes de la suite, notons-la $t_n$, donnant le nombre de structures connectives sur $X_n$ \`a isomorphisme connectif pr\`es : $t_1=1$, $t_2= 2$, $t_3= 6$, $t_4= 47$, $t_5= 3095$, $t_6= 26015236$.

\section {Representation par entrelacs}\label{Representation par entrelacs}

\begin{figure} 
\begin{center}
\includegraphics [scale=0.2]{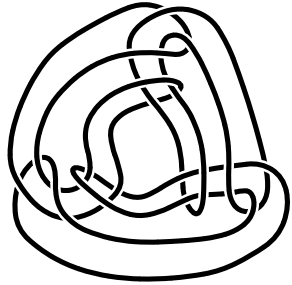}
\caption{Entrelacs de Brunn \`a trois composantes}
\label{brunn3}
\end{center}
\end{figure}

On ne consid\`ere dans cet article que les entrelacs «~apprivois\'es~» (\emph{tame links}). Un entrelacs $\tilde{E}$ \`a $n≥1$ composantes est la classe d'\'equivalence par isotopie ambiante d'un plongement $E$ dans $\mathbf{R}^3$ (ou dans la sph\`ere \`a trois dimensions $S^3$) de l'union finie disjointe de $n$ copies, num\'erot\'ees de $1$ \`a $n$, du cercle $S^1$. C'est aussi la classe d'\'equivalence de $E$ \emph{modulo} les hom\'eomorphismes de l'espace entier conservant l'orientation, que nous appellerons simplement les \emph{hom\'eomorphismes ambiants}. Un entrelacs est donc d\'efini par la donn\'ee simultann\'ee de $n$ plongements disjoints du cercle qui constituent les composantes (du repr\'esentant $E$) de l'entrelacs $\tilde{E}$, que nous noterons $E_1, E_2, \dots, E_n$. Nous identifierons  $E$ au $n$-uplet $(E_1, E_2, \dots, E_n)$. On note $\mathcal{E}_n$ l'ensemble des entrelacs \`a $n$ composantes, et $\mathcal{E}=\bigcup_{n≥1}\mathcal{E}_n$ l'ensemble des entrelacs. Etant donn\'e une partie non vide $I$ de $X_n=\{1,\dots, n\}$, on note $\widetilde{E_I}$ l'entrelacs \`a $card(I)$ composantes d\'efini par $E_I=(E_i)_{i\in I}$.

Etant donn\'e un plongement particulier $E$ d\'efinissant un entrelacs $\tilde{E}$ et $\phi$ un hom\'eomorphisme ambiant, nous noterons simplement $\phi(E)=(\phi(E_i))_{1≤i≤n}$ le plongement compos\'e $\phi\circ E$. Nous dirons que le plongement $E$ est inclus dans une partie $T$ de l'espace pour exprimer l'inclusion $E(\bigsqcup_i S^1)=\bigcup_i E_i(S^1)\subset T$.

Nous dirons qu'un entrelacs $\tilde{E}\in\mathcal{E}_n$ est \emph{s\'eparable} (\emph{splittable}) s'il existe une partition de $X_n$ en deux ensembles disjoints non vides $I$ et $J$, un hyperplan (ou, si l'on travaille dans $S^3$, une sph\`ere) $H$ et un hom\'eomorphisme ambiant $\phi$ de l'espace tels que $\phi(E_I)$ et $\phi(E_J)$ se trouvent de part et d'autre de $H$. Si l'on peut de la m\^eme fa\c con s\'eparer toutes les composantes de $E$, celui-ci est dit \emph{compl\'etement s\'eparable}. Un entrelacs non s\'eparable sera dit \emph{ins\'eparable}.

\begin{df} [Connectivity space associated to a link] La \emph{structure connective} d'un entrelacs \`a $n$ composantes $\widetilde{E}$ est la structure connective de l'espace connectif $(X_n,\mathcal{K}_{\tilde{E}})$ d\'efinie par
\begin{displaymath}
{\mathcal{K}_{\tilde{E}}}^*=\{I\in\mathcal{T}(X_n), \widetilde{E_I} \textrm { est ins\'eparable}\}.
\end{displaymath}
Nous dirons aussi que $(X_n,\mathcal{K}_{\tilde{E}})$ est l'\emph{espace connectif associ\'e} \`a l'entrelacs $\tilde{E}$, ou que l'entrelacs $\tilde{E}$ \emph{repr\'esente} l'espace connectif $(X_n,\mathcal{K}_{\tilde{E}})$. 
\end{df}

\begin {df} L'\emph{ordre connectif} d'un entrelacs est l'ordre de l'espace connectif associ\'e.
\end {df}

\noindent\emph{Exemples}. 1. La structure connective d'un n\oe ud quelconque, c'est-\`a-dire d'un entrelacs ayant une seule composante, est celle de l'espace connectif \`a un point.

\noindent 2. En 1892, Hermann Brunn~\cite{Brunn:1892a} a consid\'er\'e pour tout $n≥1$  l'entrelacs construit selon le m\^eme principe que l'entrelacs \`a trois composantes repr\'esent\'e sur  la figure~\ref{brunn3}. Nous appellerons cet entrelacs l'entrelacs de Brunn \`a $n$ composantes. En 1961, reprenant l'\'etude topologique des entrelacs de Brunn, Hans Debrunner~\cite{Debrunner:1961} a consid\'er\'e plus g\'en\'eralement ce qu'il a appel\'e les «~entrelacs de type brunniens~» ou plus simplement «~entrelacs brunniens~». En termes connectifs, les entrelacs brunniens sont pr\'ecis\'ement ceux dont la structure connective est celle que, pour cette raison, nous avons qualifi\'ee de  brunnienne (d'ordre 1). L'entrelacs borrom\'een (figure~\ref{borro}) est distinct de l'entrelacs de Brunn \`a trois composantes, mais c'est aussi un entrelacs brunnien. Le psychanalyste fran\c cais Jacques Lacan~\cite{Lacan:1973} a consid\'er\'e successivement ces deux entrelacs pour illustrer la borrom\'eanit\'e. Ainsi, de m\^eme que la notion de «~type brunnien~» introduite par Debrunner, l'id\'ee de borrom\'eanit\'e chez Lacan semble bien \^etre de nature connective. C'est ce qui nous a conduit, dans  l'article~\cite{Dugowson:2007c}, a poser la d\'efinition suivante.

\begin {df}  On appelle \emph{espaces connectifs lacaniens} les espaces connectifs admettant une repr\'esentation par entrelacs.
\end{df}

\begin{df} \label{collier}  On appelle \emph{collier} \`a $n$ composantes tout couple $(E,T)$, o\`u $T\subset \mathbf{R}^3$ est un tore solide et $E=(E_1,\dots,E_n)$ d\'efinit un plongement de l'entrelacs $\tilde E$ \`a l'int\'erieur de $T$ tel que
\begin{itemize}
\item $E$ n'est pas contenu dans une partie simplement connexe de $T$,
\item pour tout  $i\in\{1,\dots,n\}$, il existe une partie simplement connexe de $T$ qui contient $(E_j)_{j\ne i}$.
\end{itemize}
\end{df}

\noindent \emph{Remarque}. Il est \'equivalent de dire que $E$ constitue une partie \emph{essentielle} de $T$ mais que, pour tout  $i\in\{1,\dots,n\}$, $(E_j)_{j\ne i}$ est une partie non essentielle de $T$ (sur la notion de partie essentielle, voir~\cite{Rolfsen:1990}, p.\,110).

\begin{figure} 
\begin{center}
\includegraphics [scale=0.2]{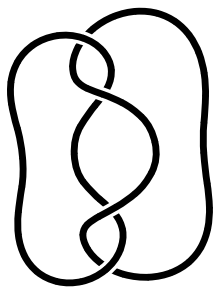}
\caption{Un collier s\'eparable}
\label{soury2}
\end{center}
\end{figure}

\noindent \emph{Exemples}. 1. L'entrelacs form\'e de deux cercles non entrelac\'es peut \^etre repr\'esent\'e par un collier (figure \ref{soury2}). Plus g\'en\'eralement, l'entrelacs compl\'etement s\'eparable \`a $n$ composantes peut \^etre repr\'esent\'e par un collier (voir par exemple, parmi les centaines d'entrelacs de toutes sortes dessin\'es par le math\'ematicien Pierre Soury celui du texte 50, page 1 de~\cite{Soury:1988.2}).

2. Les approximations finies du collier d'Antoine (\emph{Antoine's Necklace}) obtenues \`a chaque \'etape de sa construction constituent autant de colliers (finis).

\begin{figure} 
\begin{center}
\includegraphics [scale=0.2]{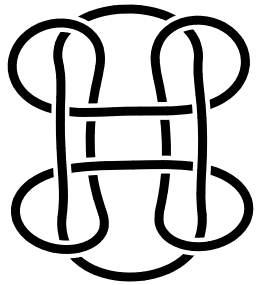}
\caption{Un diagramme du collier de Brunn \`a deux composantes}
\label{brunn2a}
\end{center}
\end{figure}

\begin{figure} 
\begin{center}
\includegraphics [scale=0.2]{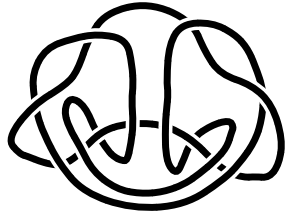}
\caption{Collier de Brunn \`a deux composantes}
\label{brunn2b}
\end{center}
\end{figure}

\begin{figure} 
\begin{center}
\includegraphics [scale=0.2]{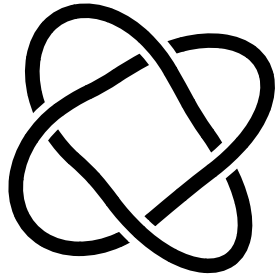}
\caption{Diagramme standard de $4_1^2$}
\label{collier2vue1}
\end{center}
\end{figure}

\begin{figure} 
\begin{center}
\includegraphics [scale=0.2]{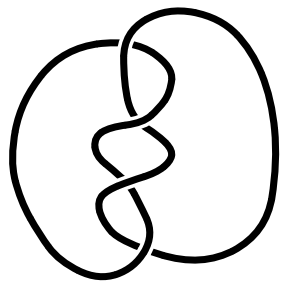}
\caption{Vue «~en collier~» de $4_1^2$}
\label{collier2vue2}
\end{center}
\end{figure}

3. L'entrelacs de Brunn \`a $n$ composantes se repr\'esente par un collier $(E,T)$ que nous appellerons le collier de Brunn. Pour $n=1$, le collier de Brunn se r\'eduit au n\oe ud trivial. Pour $n=2$, on obtient (apr\`es l'action d'un isotopie ambiante ad\'equate) l'entrelacs repr\'esent\'e sur les figures \ref{brunn2a} et \ref{brunn2b}. Remarquons que le collier de Brunn \`a deux composantes n'est pas le collier brunnien \`a deux composantes le plus simple, l'entrelacs nomm\'e $4_1^2$ dans la table de Tait \'etant lui aussi un tel collier (voir les figures \ref{collier2vue1} et \ref{collier2vue2}). 

\begin{figure} 
\begin{center}
\includegraphics [scale=0.4]{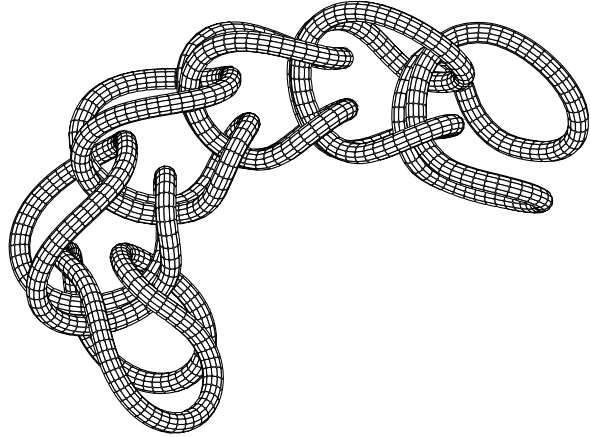}
\caption{Quelques composantes d'un collier de Brunn solide}
\label{brunnn}
\end{center}
\end{figure}

Soit maintenant $F$ le plongement dans $T$ de $n$ copies du tore solide de r\'ef\'erence $T_0=D^1\times S^1$ obtenu en rempla\c cant chaque composante $E_i(S^1)$ du collier de Brunn $E$ par un voisinage tubulaire $F_i(T_0)$ de rayon suffisament petit. Le plongement $F$ constitue ainsi un «~entrelacs brunnien de $n$ tores solides~», que nous appellerons le \emph{collier de Brunn solide} \`a $n$ composantes (la figure~\ref{brunnn} repr\'esente quelques composantes formant une partie d'un tel collier).

Dans la suite, $A$ et $B$ d\'esignant deux parties de l'espace, nous dirons que $A$ \emph{tranche} $B$ s'il existe une partie simplement connexe de $B$ contenant $B\smallsetminus A$. Il est imm\'ediat que toute partie de l'espace qui tranche le tore solide dans lequel est inscrit un collier tranche \'egalement une composante au moins de ce collier. Nous admettrons que cette propri\'et\'e continue d'\^etre v\'erifi\'ee pour les colliers de Brunn \emph{solides} : toute partie de l'espace qui tranche le tore solide dans lequel est inscrit un collier de Brunn solide tranche \'egalement une composante au moins de ce collier.

\begin{thm} L'union brunnienne d'une famille finie d'espaces connectifs finis repr\'esentables par des colliers est elle-m\^eme repr\'esentable par un collier. En particulier, les espaces brunniens it\'er\'es sont lacaniens.
\end{thm}

\noindent \emph{D\'emonstration}. Soit $(\mathbf{X}_i)_{1≤i≤p}$ une famille de $p$ espaces connectifs finis, chacun d'eux \'etant repr\'esentable par un collier $(E^i,T_i)$, avec $E^i=(E_1^i,\dots,E_{n_i}^i)$, o\`u $n_i$ est le nombre de points de l'espace $\mathbf{X}^i$. On peut chosir chaque tore solide $T_i$ de fa\c con \`a ce qu'il co\"\i ncide avec la $i$-\`eme composante d'un collier de Brunn solide \`a $p$ composantes inscrit dans un tore solide $T$. La famille $E=(E_j^i)$ o\`u $1≤i≤p$ et $1≤j≤n_i$ constitue alors l'entrelacs cherch\'e. En effet, le compl\'ementaire d'une partie simplement connexe $V$ du tore $T$  tranche n\'ecessairement l'un des $T_i$, donc l'un des $(E_j^i)$, de sorte que la famille $E$ ne peut \^etre contenue dans une partie simplement connexe de $T$. De plus, si $(E_{j_0}^{i_0})$ d\'esigne un n\oe ud quelconque composant $E$, la famille $(E_j^i)_{(i,j)\ne(i_0,j_0)}$ est constitu\'ee d'une part des entrelacs $E^i$ pour $i\ne i_0$ et d'autre part des $E_j^{i_0}$ pour $j\ne j_0$. Les entrelacs consid\'er\'es \'etant brunniens, on peut alors d\'eplacer par isotopie ambiante ces diff\'erentes composantes de fa\c con \`a mettre en \'evidence leur inclusion dans une partie simplement connexe du tore solide $T$. Ainsi, $E$ est bien un collier. Enfin, d'apr\`es ce qui pr\'ec\`ede, les seules sous-familles non vides ins\'eparables de $E$ sont les $E^i$ et $E$ elle-m\^eme, de sorte que la structure connective de $E$ est bien celle de $\biguplus_{1≤i≤p}\mathbf{X}_i$. 
La seconde affirmation en d\'ecoule imm\'ediatement, par r\'ecurrence triviale.
\begin {flushright}
$\square$
\end{flushright}
\pagebreak[3]

\begin{figure} 
\begin{center}
\includegraphics [scale=0.3]{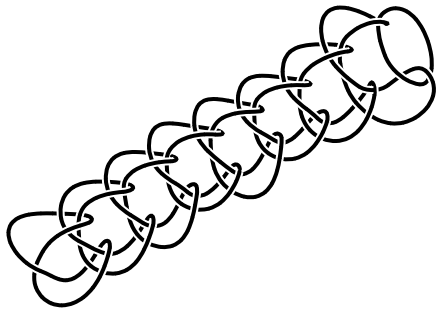}
\caption{Un entrelacs d'ordre 8}
\label{tendu}
\end{center}
\end{figure}

\noindent\emph{Exemples}. 1. Un entrelacs \`a $n$ composantes a un ordre connectif inf\'erieur ou \'egal \`a $n-1$. Voici un exemple o\`u l'ordre maximal est atteint. Soit $\mathbf{X}$ l'espace connectif de support $\{x_1,\dots, x_n\}$ dont le graphe g\'en\'erique comporte, outre les singletons, $n-1$ points g\'en\'eriques $y_1,\dots, y_{n-1}$ et dont l'ensemble des ar\^etes est  $\mathcal{A}=\{(y_i,y_{i-1}), i=1,\dots,n-1\}\cup\{(y_i,x_{i+1}), i=1,\dots,n-1\}$, o\`u on a pos\'e $y_0=x_1$. $\mathbf{X}$ est un espace brunnien d'ordre $n-1$, c'est donc un espace lacanien repr\'esentable par un entrelacs $E=(E_1, E_2,\dots, E_{n})$ pr\'esentant la propri\'et\'e suivante : si l'on coupe $E_k$, la famille $(E_1,\dots,E_{k-1})$ reste entrelac\'ee tandis que les $E_{k+1},\dots, E_n$ se s\'eparent les uns des autres. La figure~\ref{tendu} repr\'esente un tel entrelacs pour $n=9$.

\begin{figure} 
\begin{center}
\includegraphics [scale=0.2]{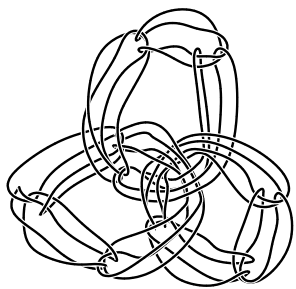}
\caption{Un «~borrom\'een de borrom\'een~»}
\label{borroborro}
\end{center}
\end{figure}

2. L'entrelacs de la figure\,\ref{borroborro} est un «~borrom\'een de borrom\'een~», d'ordre connectif \'egal \`a 2. 

3. La structure connective sur $\{a,b,c,a',b',c'\}$ qui admet pour points g\'en\'eriques $u=\{a,b\}$, $v=\{b,c\}$, $u'=\{a',b'\}$, $v'=\{b',c'\}$ et $x=\{u,v,u',v'\}$ n'est pas brunnienne, mais est facilement repr\'esentable par un collier.

\subsection{Une conjecture}\label{Une conjecture}

\emph{Conjecture}. Tout espace connectif fini est lacanien.

J'ai formul\'e cette conjecture pour la premi\`ere fois lors de  la conf\'erence~\cite{Dugowson:2007b}, puis dans l'article~\cite{Dugowson:2007c}.


\newpage

\tableofcontents

\end{document}